\documentclass[aps,pre,reprint,showpacs,groupedaddress]{revtex4-1}

\usepackage[dvips]{epsfig}
\usepackage{amsmath,amsfonts}
\usepackage{soul}
\usepackage{color}
\graphicspath{{./figures/}}

\begin{document}

\title{Excluded-volume effects in the diffusion of hard spheres}
\author{Maria Bruna}
\author{S. Jonathan Chapman}
\affiliation{University of Oxford, Mathematical Institute, 24-29 St. Giles', Oxford, OX1 3LB, United Kingdom}


\begin{abstract}
Excluded-volume effects can play an important role in determining
transport properties in diffusion of particles. Here, the diffusion of
finite-sized hard-core interacting particles in two or three
dimensions is considered 
systematically using the method of matched asymptotic expansions. The
result is a nonlinear diffusion equation for the one-particle
distribution function, with excluded-volume effects enhancing the
overall collective diffusion rate. An expression for the
effective (collective) diffusion coefficient is
obtained. Stochastic simulations of the full particle   
system are shown to compare well with the solution of this equation
for two examples. 
\end{abstract} 


\maketitle

\section{Introduction}
Recently there has been an increasing interest in
understanding  the transport of particles with size-exclusion 
\cite{sun2007toward}. 
Size exclusion is
important in many biological processes, including diffusion through ion
channels \cite{hille2001ion,gillespie2002coupling} and in chemotaxis
\cite{lushnikov2008macroscopic}, and can have
a significant  impact on the thermodynamics and kinetics of biological
processes such as association reactions at membranes \cite{kim2010crowding}. 
Finite-size effects are also important when considering the combustion
of powders \cite{saxena1990}, 
collective behavior (e.g. animal flocks or traffic movement)
\cite{camazine2003self, schadschneider2002traffic} and granular gases
\cite{barrat2005granular}.

\begin{figure}[th]
\begin{center}
\includegraphics[width=\columnwidth]{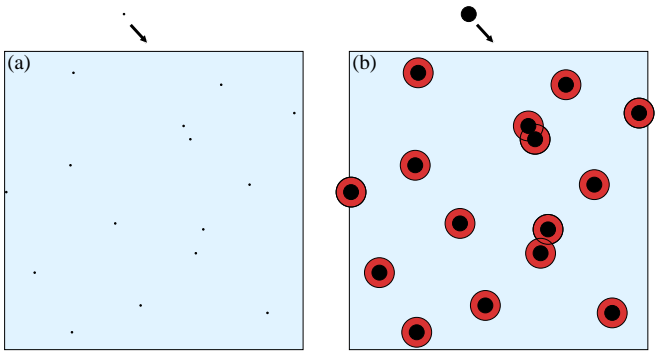}

\caption{(Color online) Excluded [red (dark grey) and black] and
  available [blue (light grey)] area in a solution of black particles for the
  placement of an additional test particle. (a) The area
  available with point particles is the whole domain. (b) The
  area available (to the \emph{center} of the test particle) with
  finite-size particles is reduced. Modified from Minton
  \cite{minton2001influence}.} 
\label{fig:steric}
\end{center}
\end{figure}

Excluded-volume or steric interactions arise from the mutual
impenetrability of finite-size particles (see Fig. \ref{fig:steric}).
For one-dimensional configurations, such as channels, the single-file
diffusion of  hard-core particles can be solved exactly by
mapping it to the classical 
diffusion of point-particles \cite{lizana2009diffusion,henle2008diffusion}. This has recently been extended to heterogeneous particles and anomalous particles \cite{flomenbom2010dynamics,flomenbom2011clustering}. 
However, the situation in higher dimensions is more challenging.

It is well-known that for finite-size particles the effective
diffusion coefficient becomes concentration-dependent. In fact we have
to distinguish between two alternative notions of diffusion
coefficient: the collective diffusion coefficient, which describes the
evolution of the total concentration, and the self-diffusion
coefficient, which describes the evolution of a single tagged particle
\cite{hanna1982self}. Here we  concentrate on  collective diffusion,
hereafter simply referred to as diffusion. 

Batchelor \cite{batchelor1976brownian} models Brownian diffusion of particles
with hydrodynamic interaction using generalized Einstein relations to
find a  concentration dependent correction to the collective diffusion
coefficient. Felderhof \cite{felderhof1978diffusion} considers the same problem
 through an analysis of the Fokker-Planck equation, and includes
both excluded volume and hydrodynamic effects. His analysis
is based on the thermodynamic limit (in which the number of particles
$N$ and the system volume $V$ tend to infinity, with the concentration
$N/V$ fixed), and is valid only for a small perturbation from the
equilibrium concentration. Similarly, the self-diffusion coefficient
to  first-order in a constant concentration is obtained from the
generalized Smoluchowski equation in
\cite{hanna1982self,ackerson1982correlations}. 

Muramatsu and Minton \cite{muramatsu1988tracer} use a simple model to calculate the
diffusion coefficient of hard spheres by estimating   the
probability that the target volume for a step in a random walk is free
of any macromolecules.
Other authors model excluded volume phenomenologically by introducing
a particle pressure, resulting in an equation of state in which the compressibility is reduced as the concentration increases \cite{degond2010congestion}. 

Another popular approach is to consider lattice models, in which a
particle can only move to a site if it is presently
unoccupied. Such an approach has been used to model diffusion of multiple species with size exclusion effects
 \cite{burger2010nonlinear, simpson2009multi} or to model the effect of crowding on diffusion-limited reaction \cite{schmit2009lattice}.

The preceding approaches are all either phenomenological in nature, 
restricted to small perturbations from a uniform concentration, or
based on the thermodynamic limit in which the number of particles
tends to infinity.
Here we consider a finite number of finite-sized particles diffusing
in a box of fixed size. We perform an asymptotic analysis
of the associated Fokker-Planck equation in the limit that the volume
fraction of particles is small. Our analysis is systematic, using the
method of matched asymptotic expansions, but is not appropriate for
concentrations close to the jamming limit.

\section{Diffusion with finite-size effects}

In order to focus on steric effects, we suppose that there are no 
electrostatic or hydrodynamic interaction forces between particles.
We work in $d$ dimensions, where $d$ is either 2 or 3.
Thus our starting point is a system of $N$ identical hard core
diffusing and interacting
spheres (or disks), each with constant diffusion coefficient $D_0$ and
diameter $K$, in a bounded domain $\Omega$ in ${\mathbb R}^d$ of typical diameter $L$. 
By nondimensionalizing length with $L$ and time with $L^2/D_0$, the
size of the domain and the diffusion coefficient may be normalized to unity, 
while the diameter of the particles becomes $\epsilon = K/L$. 
We assume that the particles occupy a small volume fraction, so that
$N \epsilon ^d \ll 1$.   
We denote the centers of the particles by ${\bf X} _i(t) \in \Omega$
at time $t\ge0$, where $1\le i \le N$ \footnote{Note that $\Omega$ is
  the space available to a particle center, which is slightly smaller
  than the container due to the finite size of the particles.}. Each
center evolves according to the stochastic differential equation (SDE)  
\begin{align}
\label{sde}
&\mathrm{d} {\bf{X}}_i \equiv  \sqrt{2}\, \mathrm{d}  {\bf B}_i
+  {\bf f}_i\, \mathrm{d} t, \qquad 1\le i \le N, 
\end{align} 
where the
${{\bf B}}_i$ are $N$ independent $d$-dimensional standard
Brownian motions and     ${\bf  f}_i$  is the external
 force on the $i$\,th particle.
In general this force
 may include both inter-particle and external interactions, such as
electromagnetic, friction, convection and potential forces, in which case
${\bf  f}_i$ depends on the positions of all the particles $\vec X
=  ({\bf  X}_1, \ldots,  {\bf  X}_N)$. 
While soft-core steric effects can also be built into ${\bf  f}_i$,
hard-core collisions can be
more easily expressed as reflective boundary conditions  on the ``collision
surfaces'' $r= ||{\bf  X}_i - {\bf  X}_j|| = \epsilon$,
with $1\le i<  j \le N$.  
Since we are focusing on hard-core particle interactions, we 
 restrict ourselves to other external forces of the form 
\begin{equation} 
\vec  F (\vec X) = [{\bf  f} ({\bf  X}_1), \dots,  {\bf f} ({\bf  X}_N)], 
\end{equation}
where ${\bf  f} : \Omega \to \mathbb R^d$ acts identically on
all $N$ particles. 
We suppose that the initial positions
${\bf  X}_i(0)$ are also random, and that they 
are independent and identically distributed.

Let  $P ({\vec x}, t)$ be the joint probability density function  
of the $N$ particles. Then, by the It\^o formula, $P  (\vec x, t)$ 
evolves according to the  linear Fokker-Planck partial
differential equation (PDE) 
\begin{subequations}
\label{fp}
\begin{align}
\label{fp_eq}
\frac{\partial P}{\partial t}&=  \vec \nabla_{\vec x} \cdot
\big[ \vec \nabla_{\vec x} \,  P -  \vec F(\vec x) \, P \big] \quad
\textrm{in}\quad \Omega_\epsilon^{N}, 
\end{align} 
where $\vec \nabla _{\vec x}$ and $\vec \nabla _{\vec x} \, \cdot$ respectively
stand for the gradient and divergence operators with respect to the
$N$-particle position vector $\vec x = ({\bf  x}_1, \dots, {\bf  x}_N) \in \Omega^{N}$. Note that because of steric effects,
\eqref{fp_eq} is not defined in  $\Omega ^N$ but in its ``hollow form''  
$\Omega_\epsilon ^N = \Omega ^N \setminus \mathcal B_\epsilon$, where
$\mathcal B_\epsilon=\{\vec x\in \Omega^N: \exists i\ne j \textrm{
  such that } ||{\bf  x}_i - {\bf  x}_j|| \le \epsilon
\}$ is the set of all illegal configurations (with at least one
overlap). On the collision surfaces $\partial \Omega_\epsilon ^N $ we
have the reflecting
boundary condition
\begin{equation} 
\label{bc} 
\big[ \vec \nabla_{\vec x} \,  P -  \vec F(\vec x) \, P \big] \cdot 
     {\vec n} = 0 \quad \text{on} \quad  \partial
     \Omega_\epsilon ^N, 
\end{equation}
\end{subequations}
where $ {\vec n}  \in \mathcal S^{dN-1}$ denotes the unit outward normal. 
Since the initial positions of the particles 
 are independent and identically distributed, the initial distribution
 function $P_0(\vec x)$ is invariant to permutations of the particle labels.
The form of (\ref{fp}) then means that  $P$ itself is invariant to
permutations  of the particle labels for all time.

Although linear, the PDE model \eqref{fp} is very high-dimensional,
and it is impractical to solve it directly. Since all the particles
are identical, we
are interested mainly in the  
 marginal distribution function of the first particle,  given by  
 $p({\bf  x}_1,t) =  \int P(\vec x,t) \  \mathrm{d}
{\bf  x}_{2} \dots \mathrm{d} {\bf  x}_{N}. $
We aim to reduce the high-dimensional PDE for $P$ to a low-dimensional
PDE for $p$ through  a systematic asymptotic expansion as $\epsilon
\rightarrow 0$.  

\subsection{Point particles}
In the particular case of point-particles ($\epsilon = 0$)  the model
reduction is straightforward. In this case the $N$ particles are
independent and the domain is $\Omega_\epsilon ^N \equiv \Omega ^N$
(no holes), which implies that the internal boundary conditions in
\eqref{bc} vanish. Therefore  $P(\vec x,t) = \prod_{i=1}^{N} p(  {\bf x}_i,t) $, and 
\begin{subequations}
\label{point}
\begin{align}
\label{point_eq}
\frac{\partial p}{\partial t}({\bf  x}_1,t) &=  {\boldsymbol \nabla}_{
  {\bf x}_1} \cdot \left[ {\boldsymbol \nabla}_{{\bf  x}_1} \,  p -  
  {\bf f}({\bf  x}_1) \, p \right] \quad  \textrm{in} \quad \Omega,\\ 
\label{point_bc} 
0&= \left[ {\boldsymbol \nabla}_{{\bf  x}_1} \,  p -  {\bf  f}(
  {\bf x}_1) \, p \right] \cdot \boldsymbol { \hat {\bf n}}_1 \quad
\textrm{on} \quad \partial \Omega, 
\end{align}
\end{subequations}
where $ \boldsymbol {\hat  {\bf n}}_1$ is the outward unit normal to
$\partial \Omega$. Note that since the particles are indistinguishable
each satisfies the same diffusion equation and boundary condition, so
that $P$ is a product of $N$ identical 1-particle distribution
functions $p$. If the particles were not identically distributed initially
then we would need a different distribution function for each one;
although these would all satisfy the same diffusion equation they
would have different initial conditions. This point will be important
when we go on to consider finite-sized particles.

\subsection{Finite-size particles}

When $\epsilon > 0$, the internal boundary conditions in \eqref{bc}
mean the particles are no longer independent. 
When we integrate \eqref{fp_eq} over ${\bf x}_2$, $\ldots$, ${\bf x}_N$
and apply the divergence theorem we end up with surface integrals 
over the collision surfaces, on which $P$ must be evaluated.
However, when the
particle volume fraction is small, the volume in  
 $\Omega_\epsilon^N$ occupied by configurations in which  three or
more particles are 
close is small [$O(\epsilon^{2d} N^2)$] compared to those in which two
particles alone are in proximity [$O(\epsilon^d N)$]. Thus the
dominant contribution to these ``collision integrals'' corresponds to
two-particle collisions.
We illustrate our approach for $N=2$; since two-particle collisions
dominate the extension to arbitrary  $N$ is  straightforward. A
similar approach is used in \cite{hanna1982self,ackerson1982correlations}.

For two particles at positions ${\bf  x}_1$ and ${\bf 
  x}_2$, Eq. \eqref{fp_eq} reads 
\begin{subequations}
\label{fp2} 
\begin{align}
\frac{\partial P}{\partial t}({\bf  x}_1, {\bf  x}_2, t)
=  &\phantom{+}{\boldsymbol \nabla}_{{\bf  x}_1} \cdot \left[ {\boldsymbol \nabla}_{{\bf x}_1}   P -  {\bf  f}({\bf  x}_1)  P
  \right]  
\nonumber\\
\label{fp2_eq} 
&+   {\boldsymbol \nabla}_{{\bf  x}_2} \cdot \left[ {\boldsymbol \nabla}_{{\bf 
      x}_2}   P -  {\bf  f}({\bf  x}_2)  P \right], 
\end{align}
for  $({\bf  x}_1, {\bf  x}_2) \in \Omega_\epsilon ^2$,
and the boundary condition \eqref{bc} reads 
\begin{equation}
\label{bc2}
\left[ {\boldsymbol \nabla}_{{\bf  x}_1}   P -  {\bf  f}(
  {\bf x}_1)  P \right] \cdot \boldsymbol {\hat {\bf n}}_1 +  \left[
  {\boldsymbol \nabla}_{{\bf  x}_2}   P -  {\bf  f}({\bf 
    x}_2)  P \right] \cdot  \boldsymbol {\hat {\bf n}}_2 = 0, 
\end{equation}
\end{subequations}
on ${\bf  x}_i \in \partial \Omega$ and $||{\bf  x}_1 -
{\bf  x}_2|| = \epsilon$. Here $ \boldsymbol {\hat {\bf n}}_i =  { {\bf n}}_i / || { {\bf n}}_i||$, where $ { {\bf n}}_i$ is the
component of the normal vector $\vec  n$ corresponding to the $i-$th
particle, $\vec n = ({\bf  n}_1, {\bf  n}_2)$. We note
that $\boldsymbol  {\hat {\bf n}}_1 = 0$ on ${\bf  x}_2 \in \partial
\Omega$, and that $\boldsymbol {\hat {\bf n}}_1 = -\boldsymbol {\hat {\bf n }}_2$ on $||
      {\bf x}_1-{\bf  x}_2|| = \epsilon$.  

We denote by $\Omega({\bf  x}_1)$ the region available to
particle 2 when particle 1 is at ${\bf  x}_1$, namely
$\Omega({\bf  x}_1) = \Omega \setminus B_\epsilon ({\bf 
  x}_1)$. Note that when the distance between ${\bf  x}_1$  and
$\partial \Omega$ is less than $\epsilon$ the volume $|\Omega(
{\bf x}_1)|$ increases. This
creates a boundary layer of width $\epsilon$ around $\partial \Omega$
where there exists a wall-particle-particle interaction (three-body
interaction). Since the dimensions of the container are much larger
than the particle diameter these interactions are higher-order and we
may safely ignore them  and take
$|\Omega({\bf  x}_1)|$ constant  
 \footnote{This would not be the case in a channel-like
  domain, for instance. In two dimensions, if $\Omega = [0, 1]\times [0, L]$ with $L
  =O(\epsilon)$, the wall-particle-particle interactions would
  become important.}.
Integrating Eq. \eqref{fp2_eq} over $\Omega ({\bf  x}_1)$ 
yields 
\begin{eqnarray} 
\lefteqn{\frac{\partial p}{\partial t}({\bf  x}_1, t) =  \phantom{-}
{\boldsymbol \nabla}_{{\bf  x}_1} \cdot \left[ {\boldsymbol \nabla}_{{\bf  x}_1} \,
  p  
-  {\bf  f}({\bf  x}_1) \, p \right] }\qquad &&\nonumber \\ 
&&\mbox{}+
 \int_{\partial B_\epsilon ({\bf  x}_1 )} [
  {\bf f}({\bf  x}_1) \, P -2{\boldsymbol \nabla}_{
  {\bf x}_1} P -{\boldsymbol \nabla}_{ {\bf  x}_2} P ]\cdot  {\boldsymbol{\hat
  {\bf n}}}_2 \, \mathrm{d} S_2 \qquad 
\label{fp2r} 
\\  &&
\mbox{ }+ \int_{\partial \Omega\cup\partial B_\epsilon ({\bf  x}_1 )}  \left[ {\boldsymbol \nabla}_{{\bf  x}_2} \,  P -
  {\bf  f}({\bf  x}_2) \, P \right] \cdot  {\boldsymbol
  {\hat {\bf  n}}}_2 \, \mathrm{d} S_2. \nonumber 
\end{eqnarray} 
The first  integral  in (\ref{fp2r}) comes from switching the
order of integration with respect to ${\bf x}_2$ and differentiation
with respect to ${\bf x}_1$ using the transport theorem; the second
 comes from using the divergence theorem on the derivatives
in ${\bf x}_2$. 
Using \eqref{bc2} and rearranging we find 
\begin{eqnarray}
\lefteqn{\frac{\partial p}{\partial t}({\bf  x}_1, t) = \phantom{-}
{\boldsymbol \nabla}_{{\bf  x}_1} \cdot \left[ {\boldsymbol \nabla}_{{\bf  x}_1} \,
  p -  {\bf  f}({\bf  x}_1) \, p \right]} &&\nonumber \\ 
\label{fp2rr}
&& \mbox{ }+ \int_{\partial B_\epsilon ({\bf  x}_1 )} \left\{- 2{\boldsymbol \nabla}_{ {\bf x}_1} P+P \left[ 
  {\bf f}({\bf  x}_1) -  {\bf  f}({\bf  x}_2)
  \right]\right\} \cdot  {\boldsymbol{\hat  {\bf n}}}_2 \, \mathrm{d} S_2.\qquad 
\end{eqnarray}
Because the pairwise particle interaction is localized near  the
collision surface $\partial B_{\epsilon} ({\bf  x}_1 )$ 
we can determine it  using the method
of matched asymptotic expansions \cite{holmes1995introduction}. 

\subsection{Matched asymptotic expansions of the density $P$}
We
suppose that when two particles are far apart ($||{\bf  x}_1\! -
{\bf  x}_2||\gg \epsilon$) they are independent, whereas when
they are close to each other ($||{\bf  x}_1 - {\bf  x}_2||
\sim \epsilon$) they are correlated. We designate these two regions of the
configuration space  $\Omega_\epsilon^2$ the outer region and inner region, respectively. 

 In the inner region, we set  $
{\bf  x}_1 = \tilde {\bf x}_1$ and ${\bf  x}_2 =
\tilde {\bf x}_1 + \epsilon \tilde {\bf x}$ and
define 
$
\tilde P (\tilde {\bf x}_1, \tilde {\bf x}, t) = P
( {\bf x}_1,  {\bf x}_2 , t)$. With this rescaling \eqref{fp2} becomes
\begin{subequations}
\label{fp_inner}
\begin{eqnarray}
\lefteqn{\epsilon^2 \frac{\partial \tilde P}{\partial t} (\tilde {{\bf 
    x}}_1, \tilde {\bf x}, t)  =   2 {\boldsymbol \nabla}_{\tilde {
    {\bf  x}}}^2 \, \tilde P  -
\epsilon ^2  {\boldsymbol \nabla} _{\tilde {\bf x}_1} \cdot \big[ 
  {\bf  f}(\tilde {\bf x}_1)  \tilde P \big ] + \epsilon^2 {\boldsymbol \nabla}_{\tilde {\bf x}_1}^2 \, \tilde P } && 
\nonumber
\\
&&\mbox{} + \epsilon  {\boldsymbol \nabla} _{\tilde {\bf x}} \cdot \big\{ \left[
  {\bf  f}(\tilde {\bf x}_1) -  {\bf 
    f}(\tilde {\bf x}_1 + \epsilon \tilde {\bf x})
  \right] \tilde P \big \}- 2 \epsilon  {\boldsymbol \nabla}_{\tilde {
    {\bf  x}}_1} \! \cdot \! {\boldsymbol \nabla}_{\tilde {\bf x}}
\tilde P ,\qquad
\label{eq_inner}
\end{eqnarray}
with 
\begin{eqnarray}
\label{bc_inner}
2 \tilde {\bf x} \cdot  {\boldsymbol \nabla}_{\tilde {\bf x}}
\tilde P &= &\epsilon \, \tilde {\bf x} \cdot \big \{
{\boldsymbol \nabla}_{\tilde {\bf x}_1} \tilde P + \left[  {\bf 
    f}(\tilde {\bf x}_1 + \epsilon \tilde {\bf x})
  - {\bf  f}(\tilde {\bf x}_1) \right] \tilde P
\big\},\qquad
\end{eqnarray}
on $|| \tilde {\bf x} || = 1$.
As noted above, we can neglect the boundary layer and hence assume that $\tilde {\bf x}_1$ is not close to $\partial\Omega$; the
region in which the
particles are close to each other and the boundary is even smaller,
and will affect only the higher-order terms.
In addition to \eqref{bc_inner} the inner solution must match with
the outer solution as $\tilde {\bf x} \rightarrow \infty$. In the outer region,
by independence, 
\begin{equation*} 
\label{indep}
P({\bf  x}_1, {\bf  x}_2, t) = q({\bf  x}_1, t)
q({\bf  x}_2, t), 
\end{equation*}
for some function $q({\bf  x}, t)$. 
Note that the invariance of $P$ with respect to a switch of particle
labels means that in the outer region both particles have the {\em
  same} distribution function $q$ (as happened in the point particle
case). 
The normalization condition on $P$ gives $q ({\bf  x}_1, t)
= p ({\bf  x}_1,t)  + \mathcal O(\epsilon^d)$. 
Expanding this outer solution in inner variables gives
\begin{eqnarray}
\label{bc_match}
P ( {{\bf 
    x}}_1, {\bf x}_2, t) &=& q(\tilde {{\bf 
    x}}_1, t) q(\tilde {\bf x}_1 + \epsilon \tilde {
  {\bf  x}})\nonumber \\
& \sim & q^2 (\tilde {\bf x}_1,t) +  \epsilon q(\tilde {
  {\bf  x}}_1) \, \tilde {\bf x} \cdot {\boldsymbol \nabla} _{ \tilde
  {\bf x}_1} q( \tilde {\bf x}_1) + \cdots.\qquad
\end{eqnarray}
\end{subequations}
Expanding $\tilde P$  in powers of $\epsilon$, and solving
(\ref{eq_inner}), (\ref{bc_inner}) with the matching condition
(\ref{bc_match}) we find 
that the solution in the inner region
is simply 
\begin{equation}
\label{sol_inner}
\tilde P (\tilde {\bf x}_1, \tilde {\bf x}, t)
\sim q^2 (\tilde {\bf x}_1,t) +  \epsilon q(\tilde {
  {\bf  x}}_1) \, \tilde {\bf x} \cdot {\boldsymbol \nabla} _{ \tilde
  {\bf x}_1} q( \tilde {\bf x}_1) + \cdots .
\end{equation}
Using this solution to evaluate the integrals in \eqref{fp2rr} we
find, to $O(\epsilon^d)$,
\begin{align}
\label{fp2_final}
\frac{\partial p}{\partial t}({\bf  x}_1, t) = \ &{\boldsymbol \nabla}_{
  {\bf  x}_1} \cdot \left [ {\boldsymbol \nabla}_{{\bf  x}_1}\! \left (  p
  +\alpha_d  \epsilon ^d p^2 \right )-  {\bf  f}(
  {\bf  x}_1) \, p \right ], 
\end{align}
where $\alpha_2 = \pi/2$ and $\alpha_3=2\pi/3$.
The extension from two particles to $N$ particles  is
straightforward up to $O(\epsilon^d)$, since at this order only pairwise
interactions need to be considered. Particle 1 has $(N-1)$ inner regions, one with each of the
remaining particles. A similar procedure shows that  
the marginal distribution function satisfies
\begin{align}
\label{fpN_final}
\frac{\partial p}{\partial t}({\bf  x}_1, t) &= {\boldsymbol \nabla}_{
  {\bf  x}_1} \!\cdot  \! \big \{ {\boldsymbol \nabla}_{{\bf  x}_1}\! \left [
  p +  \alpha_d (N \! - \! 1) \epsilon ^d p^2 \right ] -  {\bf 
  f}({\bf  x}_1) \, p \big \}.
\end{align}
Equation \eqref{fpN_final} describes the probability distribution
function for finding the 
first particle at position $x_1$ at time $t$. Since the system is
invariant to permutations of the particle labels, the marginal
distribution function of any other particle is the same.
Thus the probability  distribution
function for finding {\em any} particle at position $x_1$
at time $t$ is simply $Np$.

\subsection{Interpretation}
We see  from \eqref{fpN_final} that steric interactions lead to a
concentration-dependent 
diffusion coefficient, with the additional term 
proportional to the excluded volume. Equation \eqref{fpN_final} is
consistent with that derived by 
Felderhof \cite{felderhof1978diffusion}, but extends it to situations
in which $p$ is not close to uniform. 
We emphasize that \eqref{fpN_final} is valid for any $N$. However, for
large $N$ such that  $N-1 \approx N$ 
we can introduce  the volume concentration  $c =
\pi N  \epsilon^d p/2d$ and rewrite 
\eqref{fpN_final} as \footnote{Of course, if we define  $c =
\pi (N-1)  \epsilon^d p/2d$ then (\ref{new1}) is valid for all
$N$, but $c$ is not the volume concentration.}
\begin{equation}
\label{new1}
\frac{\partial c}{\partial t}({\bf  x}_1, t) = {\boldsymbol \nabla}_{
  {\bf  x}_1} \!\cdot  \! \big \{ D(c)\, {\boldsymbol \nabla}_{{\bf  x}_1}\!
 c  -  {\bf 
  f}({\bf  x}_1) \, c \big \},
\end{equation}
where $D(c)$ is the concentration-dependent collective diffusion
coefficient, given by 
\begin{equation}
\label{D(c)}
D(c) = 1 + 4(d-1) c.
\end{equation} 

Note that the collective diffusion coefficient $D(c)$ is increased
relative to point particles. This is in contrast to the self-diffusion
coefficient (which may be related to the mean squared displacement of
a particular particle) which is reduced relative to point particles
\cite{hanna1982self}. This apparent contradiction may be understood as
follows: the diffusion of any particular particle is impeded by its
collisions with other particles. However, these collisions bias the
random walk towards areas of low particle density, so that the overall
spread of all particles is faster. To analyze the self-diffusion
coefficient in the current framework we would need to label a
particular particle, rather than treating all particles as identical.

Whilst the self-diffusion coefficient can be thought of a diffusion
coefficient intrinsically attached to each particle, the collective
diffusion coefficient relates the diffusive flux to the concentration
gradient
of all particles \cite{mazo2009brownian}: the distribution function
$p$ is the probability of finding 
{\em any} particle at a given position, rather than the probability of
finding a particular particle there. 
 Thus the collective diffusion
coefficient is not associated with an individual (tagged) particle or
 even  a representative particle. This also means that it cannot
 easily be related to the mean-squared displacement of particles.
This distinction has important consequences when upscaling from
individual to collective behavior. 

In  \eqref{fpN_final} we have only included the leading-order
nonlinear term due to steric effects. There will be correction terms
of $O(\epsilon^{d+1}N)$ due to higher-order terms in the two-particle inner
solution (\ref{sol_inner}), as well as new inner regions where three particles [$O(\epsilon^{2d} N^2)$],
or two particles and the boundary [$O(\epsilon^{d+1}N)$], are
close. The most important of these corrections is that due to
interactions between three (or more) particles. 
Because our asymptotic expansion is systematic, these correction terms could in
principle be calculated. 

\section{Comparison with the full particle system}
In order to assess the validity of 
\eqref{fpN_final} we compare its solution $p({\bf  x}_1, t)$
(obtained by a simple finite difference method) with Monte Carlo
(MC) simulations of the $2N$-coupled SDE
\eqref{sde} in two dimensions. The particle-particle (and
particle-wall) overlaps are 
treated as in
\cite{strating1999brownian}. To test the importance of
steric interactions, we also compare with the corresponding  solutions
with $\epsilon=0$.

\begin{figure}[bh]
\begin{center}
\includegraphics[width=\columnwidth]{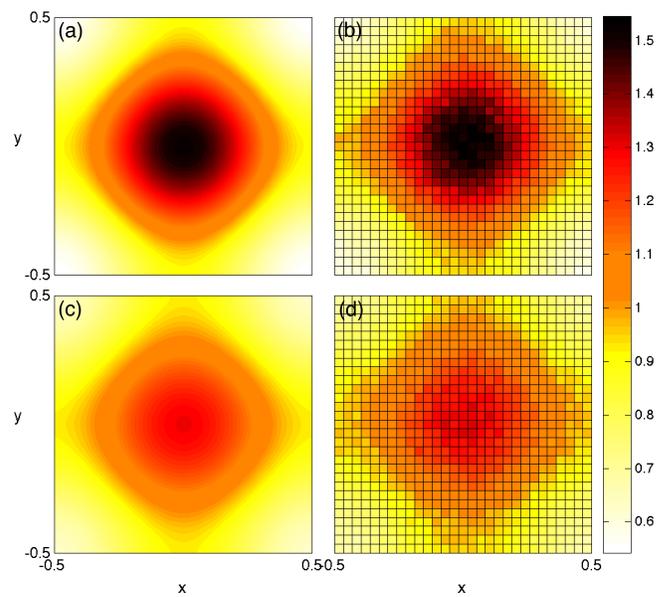}
\caption{(Color online) Marginal distribution function $p({\bf  x}_1, t)$
  at time $t=0.05$ with normally distributed initial data 
   and $N=400$. 
(a) Solution  $p({\bf  x}_1, t)$ of \eqref{point} for point
  particles ($\epsilon=0$). (b) Histogram for  $\epsilon=0$.
 (c) Solution  $p({\bf  x}_1, t)$ of \eqref{fpN_final} for
  finite-sized particles ($\epsilon=0.01$).
 (d) Histogram for $\epsilon=0.01$.
 Histograms computed from $10^4$ realizations of \eqref{sde} with $\Delta t = 10^{-5}$. All
four plots have the same color bar.} 
\label{fig:example1}
\end{center}
\end{figure}

In Fig. \ref{fig:example1} we show the results of a  time-dependent
simulation with $ {\bf  f} \equiv 0$, $\Omega = [-\tfrac{1}{2},\tfrac{1}{2}]^2$, $\epsilon
= 0.01$, $N=400$, for
which the initial 
distribution is a Gaussian of zero mean  and standard deviation 0.09
(normalized so that its integral over 
$\Omega$ is one); the figures correspond to time $t=0.05$. The simulation
time-step $\Delta t$ is chosen such that almost no collisions are
missed. 
The theoretical predictions for both point and
finite-size particles compare very well with their simulation
counterparts, while   
steric effects are clearly appreciable even though the volume
fraction of particles is only $0.0314$. However, note that while the average concentration is low, the local concentration is considerably high at the origin: $c=0.617$ at time $t=0$ and $c=0.0479$ at time $t=0.05$. The initial profile, in which particles are concentrated in the center, spreads faster when steric
effects are included [Fig. \ref{fig:example1}(c)] than when they are not
[Fig. \ref{fig:example1}(a)], indicating that the overall collective diffusion is enhanced.

\begin{figure}[hbt]
\begin{center}
\includegraphics[width=\columnwidth]{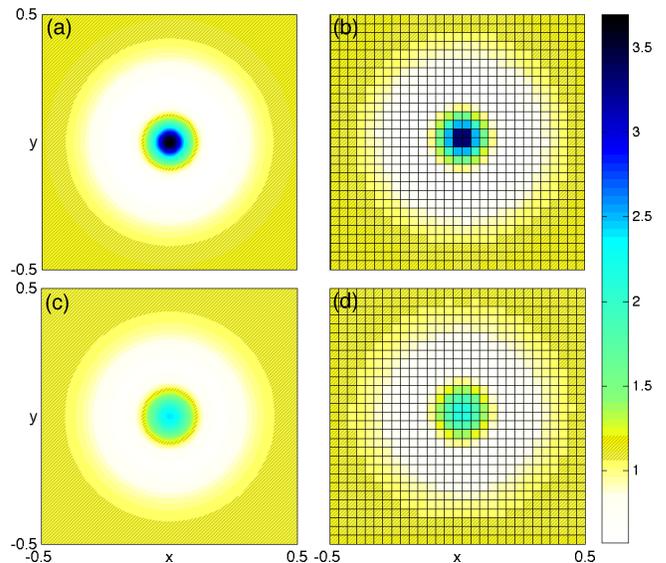}
\caption{(Color online) Stationary marginal distribution function
  $p_s({\bf x}_1)$ under the 
  external potential $V$ for point particles and finite-size
  particles, with $N=1000$.  
(a) Point particles,  $p_s \propto e^{-V}$. (b) Histogram for $\epsilon =0$. (c) Finite-size particles $p_s$ from
  \eqref{stationary} ($\epsilon = 0.01$).
 (d) Histogram for $\epsilon = 0.01$.
 Histograms computed with $10^9$ steps of the MH algorithm.  
 All four plots have the same color bar.} 
\label{fig:example2}
\end{center}
\end{figure}

When the force field ${\bf  f}$ is the gradient of a potential,
${\bf  f} ({\bf  x}_1)  = - {\boldsymbol \nabla}_{{\bf  
    x}_1} V({\bf  x}_1)$,  we may write \eqref{fpN_final} as 
\begin{equation}
\label{gradientflow}
 \frac{\partial
  p}{\partial t}  =  
{\boldsymbol \nabla}_{{\bf  x}_1} \cdot (p {\bf  u}),
\end{equation}
with $
{\bf  u} =  {\boldsymbol \nabla} _{{\bf  x}_1} \left [ \ln p + 2
        \alpha_d (N-1)\epsilon^d p + V({\bf  x}_1)\right]$. 
Equation \eqref{gradientflow} has an associated free energy \cite{carrillo2003kinetic}
\[
F(p) =\! \int_{\Omega} p \ln p + \alpha_d (N\!-\!1)\epsilon^d p^2  \, \mathrm{d} {\bf x}_1 + \int_{\Omega} \! V({\bf x}_1) \, p  \, \mathrm{d} {\bf x}_1,
\]
where the first integral corresponds to the internal energy and the second integral is the potential energy. Note that excluded-volume effects increase the internal energy of the system. 
The stationary distribution, which we denote $p_s({\bf  x}_1)$, is obtained by minimizing the free energy or by solving 
\begin{equation}
\label{stationary}
\ln p_s({\bf  x}_1) + 2\alpha_d (N-1)\epsilon^d p_s({\bf  x}_1) + V({\bf  x}_1) = C,
\end{equation} 
with the constant $C$ determined by the normalization condition on
$p_s$. 
For our second example we 
consider the volcano-shaped potential  
$V({\bf  x}_1) = -4.77 \, e^{-100||{\bf 
    x}_1||^2}+ 3.58\,  e^{-50 ||{\bf  x}_1||^2}$
 and we
compare the 
stationary distribution $p_s$ predicted by \eqref{stationary} with simulations
using the  Metropolis-Hastings (MH) algorithm
\cite{chib1995understanding}. Figure \ref{fig:example2} shows the model and simulation results with $N=1000$ and $\Omega$ and $\epsilon$ as in Fig. \ref{fig:example1} for both  point and finite-size particles (with volume fraction 0.079 and volume concentration $c=0.189$ at the origin). In this case there is competition between the potential well and steric repulsion: the particle density inside the well is reduced for finite-size particles. Again, the agreement between the model \eqref{stationary} and the stochastic simulations is excellent.

\section{Discussion}
We have derived systematically a nonlinear diffusion equation which
describes steric interactions in the limit of small but finite
particle volume fraction. Our method justifies for example the ansatz
made in \cite{calvez2006volume}  to account for the finite size of the
cells in the Keller-Segel model and prevent aggregation, and unlike \cite{felderhof1978diffusion,
  beenakker1983self}  does not rely on a closure assumption.

The equation we have derived is for the one-particle distribution
function, which measures the probability of finding {\em any} particle
at a given position; the particles we consider are identical and
indistinguishable. This means that we are examining {\em  collective
diffusion}, and we find that this is enhanced by the finite size of the
particles.
We have not considered the self-diffusion of a particular (tagged)
particle, which can be related to an individual particle's mean-square
displacement.  To analyze the self-diffusion
coefficient in the current framework we would need to label a
particular particle, rather than treating all particles as identical;
we intend to do this in a future work in which we consider multiple
particle populations.  

We note that for point particles in one dimension 
(where particles must move in single file and are not allowed to pass)
\cite{flomenbom2008single} has 
observed density dependence in the self-diffusion (mean square
displacement) of a tagged particle. Their interpretation is that the
expansion of the whole system from dense to dilute environments
quickens the self-diffusion of any tagged particle.
This is a different effect to the one we have observed, since, as
mentioned in the introduction, the {\em collective} diffusion of 
point particles in one dimension is linear, with a diffusion
coefficient independent of the
density. 
In  two dimensions self-diffusion is less sensitive to the particle
density, since, informally, there is much more space for particles to
pass each other.

The method we have developed was implemented here in its simplest setting
(hard-core identical spherical particles with an external potential)
but it can be extended in many directions. Particles of different size (\emph{cf.} \cite{bruna2011excluded}) 
or shape  can easily be incorporated, while the
hard-core interaction between particles can be replaced by any
short-range  soft-core
interaction.

On the other hand, incorporating long range effects such as chemotaxis
or electrostatic interactions is more challenging; in such cases a
system size expansion is likely to be needed in addition to a small
particle expansion.

\acknowledgments

This publication was based on work supported in part by Award No. KUK-C1-013-04, made by King Abdullah University of Science and Technology (KAUST). M.B. acknowledges financial support from EPSRC and the program ``Partial Differential Equations in Kinetic Theories'' at the Isaac Newton Institute for Mathematical Sciences, Cambridge, UK. The authors also thank Jos\'e Antonio Carrillo for valuable discussions.

\end{document}